\documentclass[twocolumn]{autart}    

\usepackage{amsfonts, amssymb, psfrag, color, multirow}
\usepackage{graphicx, subfigure}     

\newtheorem{lemma}{\textbf{Lemma}}
\newtheorem{theorem}{\textbf{Theorem}}

\def\QEDMark{\square}

\begin{document}

\begin{frontmatter}

\title{A Stochastic Kaczmarz Algorithm for Network Tomography\thanksref{footnoteinfo}}

\thanks[footnoteinfo]{The work of V.~Borkar is supported in part by a J. C. Bose Fellowship from the Government of India, grant 11IRCCSG014 from IRCC, IIT Bombay, and a grant from the Dept. of Science and Technology for `Distributed Computation for Optimization over Large Networks and High Dimensional Data Analysis'. D. Manjunath is also affiliated with the Bharti Centre for Communication, IIT Bombay and is also supported in part by the above DST grant. This paper was not presented at any IFAC meeting. Corresponding author G.~Thoppe. Tel. +91-22-22782930.\\ The authors thank the referees for their constructive comments that have helped us to improve the article significantly.}

\author[TIFR]{Gugan Thoppe}\ead{gugan@tcs.tifr.res.in},
\author[IITB]{Vivek Borkar}\ead{borkar.vs@gmail.com},
\author[IITB]{D. Manjunath}\ead{dmanju@ee.iitb.ac.in}

\address[TIFR]{School of Technology and Computer Science, Tata Institute of Fundamental Research, Mumbai 400005, India}
\address[IITB]{Department of Electrical Engineering, Indian Institute of Technology, Powai, Mumbai 400076, India}


\begin{abstract}
We develop a stochastic approximation version of the classical Kaczmarz algorithm that is incremental in nature and takes as input noisy real time data.  Our analysis shows that with probability one it mimics the behavior of the original scheme: starting from the same initial point, our  algorithm and the corresponding deterministic Kaczmarz algorithm converge to precisely the same point. The motivation for this work comes from network tomography where network parameters are to be estimated based upon end--to--end measurements. Numerical examples via Matlab based simulations demonstrate the efficacy of the algorithm.
\end{abstract}

\end{frontmatter}

\section{Introduction}
\label{sec:intro}
\subsection{Kaczmarz algorithm}

Kaczmarz algorithm \cite{Kaczmarz37} is a successive projection based iterative scheme for solving ill posed linear systems of equations. Since its introduction, its convergence properties have been extensively analyzed \cite{Galantai04} and it has found diverse applications in areas ranging from tomography \cite{Popa04}, synchronization in sensor networks \cite{Freris12}, to learning and adaptive control \cite{Astrom83,Parks92,Richalet78}. The original algorithm is deterministic, but some applications, notably network tomography which we describe later, call for a stochastic version. In this article, we introduce and analyze a stochastic approximation version based on the  Robbins-Monro paradigm \cite{Robbins51} that has become a standard workhorse of signal processing and learning control \cite{Benveniste90}, \cite{Kushner03}. We use the `o.d.e.' approach \cite{Derevitskii74,Ljung77} to analyze the scheme and argue that it has the same asymptotic behavior as the original deterministic scheme `almost surely'. While we apply our results to network tomography in this article, we believe that this analysis will be of use in other areas mentioned above. In particular, networked control is one potential application area.

A significant development in this line of research has been a randomized Kaczmarz scheme having provable strong convergence properties \cite{Strohmer09}, \cite{Leventhal10}, with recent modifications to further improve performance by weighted sampling \cite{Freris12}, \cite{Zouzias12}. The important difference between these works and ours is as follows. For them, the randomization is over the choice of rows, which is a part of algorithm design and can be chosen at will. In our case, however, a part of the randomness is due to noise and not under our control, as also in the choice of rows which a priori we allow to be uncontrolled.

\subsection{Network tomography}

Network tomography is inference of spatially localized network behavior using only measurements of end-to-end aggregates. Recent work can be classified into traffic volume and link delay tomography. A basic paradigm in both these is to infer the statistics of the random vector $X$ from an ill posed measurement model $Y = AX,$ where the matrix $A$ is assumed to be known a priori. See \cite{Coates02b,Castro04,Lawrence07} for excellent surveys.

In the transportation literature, the aim is to estimate the traffic volume on the end-to-end routes assuming access to only traffic volumes on a subset of links \cite{Maher83,Bell91,Sherali94,Vardi96}. An excellent survey is given in \cite{Abramsson98}. An analogous problem has been addressed in packet networking \cite{Feldmann01,Zhang03a,Zhang03b}. In all of these works, one sample of $Y$ is assumed available and $X$ is estimated by a suitable regularization.

Link delay tomography deals with estimation of link delay statistics from path delay measurements. Here the network is usually assumed to be in the form of a tree. Multicast probe packets, real or emulated, are sent from the root node to the leaves. For each probe packet, a set of delay measurements for paths from the root node to the leaf nodes is collected. These delays are correlated and this correlation is exploited to estimate the link delay statistics \cite{Caceres99,Adams00,Duffield01,LoPresti02}. Using $T$ independent samples of the path delay vector, an expectation maximization based algorithm is derived to obtain the maximum likelihood estimates of the link parameters. There is also work on estimating link level loss statistics \cite{Caceres99}, link level bandwidths \cite{Downey99}, link-level cross traffic \cite{Prasad03} and network  topology \cite{Duffield02} using end-to-end measurements.

\subsection{Summary of our work}

We show that, starting from the same initial point, the stochastic approximation variant of the Kaczmarz (SAK) algorithm and the deterministic Kaczmarz algorithm converge to the same point. Using this, we develop a novel online algorithm for estimation of the means (more generally, moments and cross-moments) of the elements of the vector $X$ from a sequence of measurements of the elements of the vector $Y = AX.$ Our scheme can be used for both traffic volume and link delay tomography. An important advantage of our scheme is that it is real time---taking observed data as inputs as they arrive and making incremental adaptation. Also, unlike previous approaches, our scheme allows for elements of $X$ to be correlated. While our analysis is under the simplifying statistical assumption that the samples are IID, we point out later in Section 5 that these can be relaxed considerably. For link delay tomography, our algorithm does away with the need for multicast probe packet measurements and can be used even for networks with topologies other than tree.


\section{Model and problem description} \label{sec:model}

\subsection{Basic notation}

For $n \in \mathbb{N},$  $[n] := \{1, \ldots, n\}.$ For vectors, we use $||\cdot||$ to denote their Euclidean norm and $\langle\cdot, \cdot \rangle$ for inner product. For a matrix $A,$ $a_i$ denotes its $i$-th row, $a_{ij}$ its $(i,j)$-th entry, $\mathcal{R}_A$ its row space and $A'$ its transpose. We use $\dot{x}(t)$ to denote the derivative of the map $x$ with respect to $t.$

Let $X \equiv (X(1), \ldots, X(N))^{\prime}$ denote the random vector with finite variance whose statistics we wish to estimate. Let $A \in\mathbb{R}^{m \times N}, m < N,$ be an a priori known matrix with full row rank and let
\begin{equation} \label{eq:LinearRel}
{
Y \equiv (Y(1), \ldots, Y(m))^{\prime} = AX + W,
}
\end{equation}
where $W$ is a zero mean, bounded variance random variable denoting noise in the measurement.
Let $Z$ be a random variable taking values in $[m]$ such that, $\forall i \in [m],$ $\Pr\{Z = i\} =: \lambda_i > 0.$ Let $\{X_k\}, \{Z_k\}, \{W_k\}, k \geq 1,$ be  IID copies of $X, Z, W$, that are jointly independent, and $Y_k := AX_k + W_k$. (The IID assumption is purely for simplicity of analysis. We point out later that these results extend to much more general situations.) We assume that at each time step $k,$ we know only the value of $Z_{k + 1}$ and the $Z_{k + 1}$-th component of $Y_{k+1},$ i.e. $Y_{k + 1}(Z_{k + 1}) =: \mathcal{Y}_{k + 1}.$

Our objective is to develop a real-time algorithm, with provable convergence properties, to estimate the moments and cross moments of the random vector $X.$

\section{Preliminaries}
\label{sec:preliminaries}

\subsection{Stochastic approximation algorithms}
\label{subsec:StocApprox}

The archetypical stochastic approximation algorithm is
\begin{equation}\label{eq:BasicStocApproxAlgo}
{
x_{k + 1} = x_k + \eta_k[h(x_k) + \xi_{k + 1}],
}
\end{equation}
where $h: \mathbb{R}^n \rightarrow \mathbb{R}^n$ is Lipschitz, $\{\eta_k\}_{k \geq 0}$ is a  positive stepsize sequence satisfying   $\sum_{k \geq 0}\eta_{k}=\infty$ and  $\sum_{k \geq   0}(\eta_{k})^{2}<\infty,$ and  $\xi_{k + 1}$ represents noise. As $\eta_k \rightarrow 0,$ (\ref{eq:BasicStocApproxAlgo}) can be viewed as a noisy discretization of the o.d.e.
\begin{equation} \label{eq:limitingODE}
\dot{x}(t) = h(x(t)).
\end{equation}
This is the `o.d.e.\ approach' \cite{Derevitskii74,Ljung77}. More specifically, suppose that the following assumptions hold.

(A1) $\{\xi_{k}\}$ is a square-integrable martingale difference sequence w.r.t.\ the $\sigma-$fields $\{\mathcal{F}_{k}\},$  $\mathcal{F}_{k}:=\sigma(x_{0},\xi_{1}, \ldots, \xi_{k}),$ satisfying $E[||\xi_{k+1}||^{2}|\mathcal{F}_{k}]\leq L(1+||x_{k}||^{2})$ a.s.\ for some $L >0.$

(A2) $\forall u,  h_\infty(u) := \lim_{c \uparrow \infty}h(cu)/c $ exists ($h_{\infty}$ will be necessarily Lipschitz) and the o.d.e. $\dot{x}(t)=h_{\infty}(x(t))$ has origin as its globally asymptotically stable equilibrium.

(A3) $H := \{x \in \mathbb{R}^n : h(x) = 0\} \neq \emptyset.$ Also, $\exists$ a continuously differentiable Lyapunov function $\mathcal{L} : \mathbb{R}^{n} \rightarrow \mathbb{R}$ such that $\langle\nabla\mathcal{L}(x), h(x)\rangle < 0$ for $x \notin H.$

Then, as in Chapters 2,3 of \cite{Borkar08}, we have:
\begin{lemma}
  \label{lem:convToIsolEq}
  The iterates $\{x_k\}$ of (\ref{eq:BasicStocApproxAlgo}) a.s. converge to $H.$
\end{lemma}

\subsection{Kaczmarz algorithm}

Consider the inverse problem of finding a fixed $v^* \in \mathbb{R}^{N}$ from $Av^*,$ where $A$ is as defined in Section~\ref{sec:model}. W.l.o.g., let rows of $A$ be of unit norm. Given an approximation $x_0$ of $v^*,$ a natural optimization problem to consider is
\begin{equation} \label{eq:optProb}
  \min_{u \in \mathbb{R}^N}\left\Vert u - x_{0}\right\Vert, \mbox{ subject to } Au = Av^*.
\end{equation}
Elementary calculation shows that its solution is
\begin{equation}
  \label{eq:sol}
  x^{*} = x_{0} + A^{\prime}(AA^{\prime})^{-1}(Av^* - Ax_{0}).
\end{equation}

Clearly, $x^* \in \mathcal{A}^{0} := x_{0} + \mathcal{R}_A.$ As $A$ has full row rank, $x^*$ is the only point in $\mathcal{A}^0$ that satisfies $Au = Av^*.$ The Kaczmarz algorithm uses this fact to solve (\ref{eq:optProb}). With prescribed initial point $x_{0},$ stepsize $\kappa,$ and $r_k \equiv (k \mbox{ mod } m) + 1,$ its update rule is given by
\begin{equation}
  \label{eq:implRule}
  x_{k+1} = x_{k} + \kappa [\langle a_{r_k} , v^* \rangle - \langle a_{r_k}, x_{k}\rangle]a_{r_k}.
\end{equation}
\begin{theorem}
  \label{thm:KaczConv}
  \cite{Chong01}
  If $0 < \kappa < 2,$ then $x_{k}\rightarrow x^{*}$ as $k \rightarrow \infty$.
\end{theorem}

Let $\mathcal{A}^{*} := v^* + \mathcal{R}_A.$ Since $\mathcal{A}^{0},$ $\mathcal{A}^{*}$ are translations of $\mathcal{R}_A,$ $\mbox{dist}(x_{0},\mathcal{A}^{*}) = \mbox{dist}(\mathcal{A}^{0},v^*).$  As $A(x^{*} - v^*) = 0,$ $(x^{*} - v^{*}) \perp \mathcal{R}_A.$ Thus, $(x^{*} - v^{*}) \perp \mathcal{A}^{0}, \mathcal{A}^{*}.$ Hence, $||v^* - x^{*}|| = \mbox{dist}(\mathcal{A}^{0},v^*) = \mbox{dist}(x_{0},\mathcal{A}^{*}).$ Thus we have:

\begin{lemma} \label{lem:QualityInitPoint}
  For any $\delta>0,$ $\left\Vert x^{*}-v^*\right\Vert < \delta$ if and only if $\mbox{dist}(x_{0},\mathcal{A}^{*})<\delta.$
\end{lemma}

\section{The SAK Algorithm} \label{sec:SAKAlg}
We develop here a SAK algorithm to estimate $\mathbb{E}X$ for the model of Section \ref{sec:model}. Let $x_0,$ an approximation to $\mathbb{E}X,$ be given. Observe from (\ref{eq:LinearRel}) that
\begin{equation}
  \mathbb{E}Y=A\mathbb{E}X.\label{eq:delayMatrixFormWithExpectation}
\end{equation}
By rescaling equations, we assume w.l.o.g. that the rows of $A$ are of unit norm. This saves some notation without affecting the analysis. $\mathbb{E}Y$ not being known exactly, one may estimate it off-line and use the classical Kaczmarz to determine $\mathbb{E}X.$ From (\ref{eq:sol}), note that the classical Kaczmarz would have converged to
\begin{equation}
\label{eq:clSol}
x^* = x_0 + A' (AA')^{-1}(\mathbb{E}(Y) - Ax_0).
\end{equation}
As against this off-line scheme, a better alternative is to use an on-line algorithm. Using the notations and assumptions of Section~\ref{sec:model}, a SAK algorithm to estimate $\mathbb{E}X,$ based on (\ref{eq:implRule}),  is:
\begin{equation} \label{eq:SAK}
  x_{k + 1} = x_k + \eta_k [\mathcal{Y}_{k+1} - \langle a_{Z_{k + 1}}, x_k \rangle] a_{Z_{k + 1}},
\end{equation}
where $\{\eta_k\}$ is as defined below (\ref{eq:BasicStocApproxAlgo}). Note in (\ref{eq:SAK}) the noisy measurements $\{\mathcal{Y}_{k}\}$ of the elements of $\mathbb{E}Y$ and the real time estimates $\{x_k\}$ of $\mathbb{E}X.$

We now analyze its behaviour. Clearly, the iterates $\{x_k\}$ of (\ref{eq:SAK}) always remain confined to $\mathcal{A}^0,$ the affine space defined below (\ref{eq:sol}). Since $A$ has full row rank, for each $k \geq 0,$ there exists unique $\alpha_k \in \mathbb{R}^m$ such that
\begin{equation} \label{eq:alphaRelx}
  x_k = x_{0} + A^\prime\alpha_k.
\end{equation}
Thus one can equivalently analyze the algorithm
\begin{equation} \label{eq:SAKEquiv}
  \alpha_{k + 1} = \alpha_{k} + \eta_k [\tilde{Y}_{k + 1} - e_{Z_{k + 1}}A(x_0 + A^\prime \alpha_k)],
\end{equation}
where $\alpha_0 = 0,$ $e_{Z_{k + 1}}$ is the ${m \times m}$ matrix with $1$ in its $Z_{k + 1}$-th diagonal position and zero elsewhere and $\tilde{Y}_{k + 1}$ is the $m-$dimensional vector with its $Z_{k + 1}$-th position occupied by $\mathcal{Y}_{k + 1}$ and zero elsewhere.

Let $\gamma_{k + 1} = [\tilde{Y}_{k + 1} - e_{Z_{k + 1}}A(x_0 + A^\prime \alpha_k)].$ Defining $  \xi_{k + 1} = \gamma_{k + 1} - \Lambda \left(\mathbb{E}Y -  A(x_0 + A^\prime \alpha_k)\right),$ where $\Lambda = \mbox{diag}(\lambda_1, \ldots, \lambda_m),$ note that (\ref{eq:SAKEquiv}) can be rewritten as
\begin{equation}
\label{eq:SAKStdForm}
  \alpha_{k + 1} = \alpha_k + \eta_k\left[\Lambda \left(\mathbb{E}Y - A(x_0 + A^\prime \alpha_k)\right) + \xi_{k + 1}\right].
\end{equation}

If $h(u) := \Lambda \left(\mathbb{E}Y - A(x_0 + A^\prime u)\right),$ then, clearly, (\ref{eq:SAKStdForm}) is in the form given in (\ref{eq:BasicStocApproxAlgo}). Its limiting o.d.e. is thus
\begin{equation} \label{eq:EquivLimODE}
  \dot{\alpha}(t) = \Lambda \left(\mathbb{E}Y - A(x_0 + A^\prime \alpha(t))\right).
\end{equation}
\begin{theorem} \label{thm:StocConv}
  $\alpha_k \stackrel{k \uparrow \infty}{\rightarrow} \alpha^* := (AA^\prime)^{-1}(\mathbb{E}Y - Ax_0).$
\end{theorem}

\textbf{Proof.} For each $k \geq 0,$ let $\mathcal{F}_{k} := \sigma(\alpha_0, \xi_{1}, \ldots, \xi_{k}).$ Lipschitz property of $h$ and (A1) are easily verified. If $h_{c}(u) := h(cu)/c$, then $h_{c}(u) \stackrel{c\uparrow \infty}{\rightarrow} h_{\infty}(u) := -\Lambda AA^{\prime}u,$ pointwise. Let $\mathcal{L}_\infty (u) := \|A^\prime u\|^2$. This vanishes only at the origin. Further, for any solution to the o.d.e. $\dot{\alpha}(t) = -\Lambda AA^\prime \alpha(t),$  $\dot{\mathcal{L}}_\infty(\alpha(t)) = - 2 ||\sqrt{\Lambda} A A^\prime \alpha(t)||^2 \leq 0,$ again with equality only at the origin. Thus $\mathcal{L}_\infty$ is a Lyapunov function for the o.d.e. $\dot{\alpha}(t) = - h_\infty(\alpha(t))$ with the origin as its  globally asymptotically stable equilibrium. Thus $(A2)$ holds. By Lemma~\ref{lem:convToIsolEq}, to exhibit $\alpha_k \rightarrow \alpha^*,$ it now suffices to show $(A3),$ i.e. $\alpha^*$ is the globally asymptotically stable equilibrium of the o.d.e. given in (\ref{eq:EquivLimODE}). Towards this, consider the function $\mathcal{L}(u) = ||A^\prime (u - \alpha^*)||^2.$ As $A$ has full row rank, $\mathcal{L}(u) = 0$ if and only if $u = \alpha^*.$ For any solution $\alpha(t)$ of (\ref{eq:EquivLimODE}), $\dot{\mathcal{L}}(\alpha(t)) = 2 \langle A^\prime (\alpha(t) - \alpha^*), A^\prime \dot{\alpha}(t) \rangle.$ But $\dot{\alpha}(t) = -\Lambda AA^\prime (\alpha(t) - \alpha^*).$ Thus $\dot{\mathcal{L}}(\alpha(t)) \leq 0$ with equality only when $\alpha(t) = \alpha^*.$ This shows that $\mathcal{L}$ is a Lyapunov function. Thus $\alpha^*$ is the sole globally asymptotically stable equilibrium of (\ref{eq:EquivLimODE}) as desired.  \qquad $\QEDMark$

Because of (\ref{eq:alphaRelx}), it follows that the SAK algorithm of (\ref{eq:SAK}) converges to $x^*$ of (\ref{eq:clSol}), the same point that the corresponding classical Kaczmarz converges to.

\section{Extensions} \label{sec:Extns}

\begin{enumerate}

\item We had assumed $\{Z_k\}$ to be IID. The final result, however, can be established under much more general conditions. For example, $\{Z_k\}$ can be:
 \begin{itemize}
 \item ergodic Markov, as in Part II, Chapter 1, \cite{Benveniste90}, with $\lambda_i$'s the corresponding stationary probabilities.

 \item asymptotically stationary as in Chapter 6, \cite{Kushner03},

 \item `controlled' Markov, as in Chapter 6, \cite{Borkar08}, which allows for non-stationarity under the mild restriction that the relative frequency of $Z_k = i$ remain bounded away from zero a.s.\ $\forall i$.
 \end{itemize}
 Likewise, $\{X_k\}$ can be ergodic Markov or asymptotically stationary as long as it is independent of $\{Z_k\}.$ In fact, we can allow it to be long range dependent and heavy tailed \cite{Ananth12}, which is often the case with real communications networks.

\item The above analysis was for estimation of means.  We can also extend it to cover higher moments. For simplicity, we neglect measurement noise from (\ref{eq:LinearRel}) and consider the model $Y = AX.$ Observe that, for any $q \in \mathbb{N}$ and each $i \in [m],$
\begin{equation} \label{eq:HighMomRel}
[Y(i)]^q = \sum_{\mathbf{r} \in \Delta_{N,q}} {q \choose {r_1, \ldots, r_N}} \prod_{\ell = 1}^N [a_{i\ell} X(\ell)]^{r_\ell},
\end{equation}
\ \\
\noindent where $\mathbf{r} \equiv (r_1, \ldots, r_N),$ ${q \choose {r_1, \ldots, r_N}} = \frac{q!}{r_1! \cdots r_N!}$ and $\Delta_{N,q} = \{\mathbf{r} \in \mathbb{Z}_{+}^N: \sum_\ell r_\ell = q\}.$ Let $\tilde{X}^q$ denote $\left( \prod_{\ell = 1}^{N}[X(\ell)]^{r_\ell}: \mathbf{r} \in \Delta_{N, q} \right)^\prime,$ a ${N + q - 1 \choose {N - 1}}$ dimensional vector,  and let $Y^q \equiv ([Y(1)]^q, \ldots, [Y(m)]^q)^\prime.$ Also, let $A^q$ denote the $m \times {N + q - 1 \choose {N - 1}}$ matrix, whose $(i,j)$-th entry is the coefficient associated with $j^{th}$ component of $\tilde{X}^q$ as given in (\ref{eq:HighMomRel}). The set of relations in (\ref{eq:HighMomRel}) can thus be compactly written as
\begin{equation} \label{eq:HighMomMatRel}
Y^q = A^q \tilde{X}^q.
\end{equation}
\vspace{2ex}
Note that $A^q$ is generically full row rank \footnote{Fix $K \geq q.$ Let $B = [[b_{ij}]]$ be an $m\times m$ non-singular matrix made up of the independent columns of $A,$ $B^k := [[(b_{ij})^k]]$. Clearly, columns of $B^k$ are also columns of $A^k,$ where $A^k$ is as  defined above (\ref{eq:HighMomMatRel}). Let $\mathcal{O}^k := \{B : det(B^k) \neq 0\},$ $k \geq 0.$ Each $\mathcal{O}^k$ contains the identity matrix, hence is non-empty. Since $det(B^k)$ is a multivariate polynomial, $\mathcal{O}^k$ is a non-empty Zariski open set in $R^{m^2}$, therefore open dense in the usual topology. Thus $\cap_{k = 0}^K\mathcal{O}^k$ is an open dense set of $R^{m^2}$.}. Hence, (\ref{eq:HighMomMatRel}) is of the same spirit as (\ref{eq:LinearRel}). One can thus use (\ref{eq:SAK}), after replacing samples of $Y_i$ with those of $Y_i^q$ and $A$ with $A^q,$ to estimate in real-time $\mathbb{E}\left(\prod_{j = 1}^N X_j^{r_j}\right)$ for any $\mathbf{r} \in \Delta_{N, q}.$ For desired $\mathbf{r},$ the only condition one needs to ensure is that if $r_{j_1}, \ldots, r_{j_\ell}$ are the components of $\mathbf{r}$ that are positive, then $\exists i \in [m]$ such that $a_{ij_1}, \ldots, a_{ij_\ell}$ are simultaneously nonzero. Clearly, by choosing appropriate $q$ in (\ref{eq:HighMomRel}), one can estimate the moments of any desired order.

\noindent Given finite moment estimates, one can then postulate a maximum entropy distribution. For e.g., if $E[\|X\|^2] \approx a,$ $E[\|X\|^4] \approx b$, then the maximum entropy distribution is $\rho^{-1}exp\left(-(\lambda \|x\|^2 + \mu\|x\|^4)\right)$, where $\rho$ is for normalization and $\lambda, \mu$ are chosen so as to ensure $E[\|X\|^2] = a, E[\|X\|^4] = b$.

\item We have taken the process $\{Z_k\}$ as given, i.e., not within our control. If instead
one can schedule $\{Z_k\}$, randomization policies such as \cite{Strohmer09,Leventhal10,Freris12,Zouzias12} can be used to advantage. Further performance improvements are possible by adapting additional averaging as in \cite{Polyak92}. We do not pursue this here.

\item Since (\ref{eq:SAKStdForm}) is of the form $\alpha_{k+1} = D_k\alpha_k + \eta_k(b + \xi_{k+1})$ for $D_k := I - \eta_k\Lambda AA^{\prime},$ suitable vector $b$ and martingale difference noise $\xi_k$, we can iterate this to obtain $\alpha_k = \prod_{m=0}^{k - 1}D_m\alpha_0 + \sum_{m=0}^{k-1}\eta_m\prod_{\ell = m+1}^{k-1}D_{\ell}(b + \xi_{m+1})$. Note that the matrix $\Lambda AA^{\prime}$ is similar to the positive definite matrix $\sqrt{\Lambda}AA^{\prime}\sqrt{\Lambda}$. Hence, if $\zeta > 0$ denotes the minimum eigenvalue of $\sqrt{\Lambda}AA^{\prime}\sqrt{\Lambda},$ then $\|\prod_{m=j}^kD_m\|_{\lambda} \leq \prod_{m=j}^k(1 - \eta_m\zeta) \leq e^{-\zeta(\sum_{m=j}^k\eta_m)}$, where $\| \ \cdot \ \|_{\lambda}$ denotes the weighted norm defined by $\|r||_{\lambda} := \left(\sum_i\frac{r_i^2}{\lambda_i}\right)^{\frac{1}{2}}$.
This can be used to obtain estimates  for finite time error and convergence rate. We do not pursue this here, see, however, \cite{Moustakides98}.

\end{enumerate}

\section{Experimental Results} \label{sec:experiments}

\begin{figure}
  \centering
  \includegraphics[scale = 0.55]{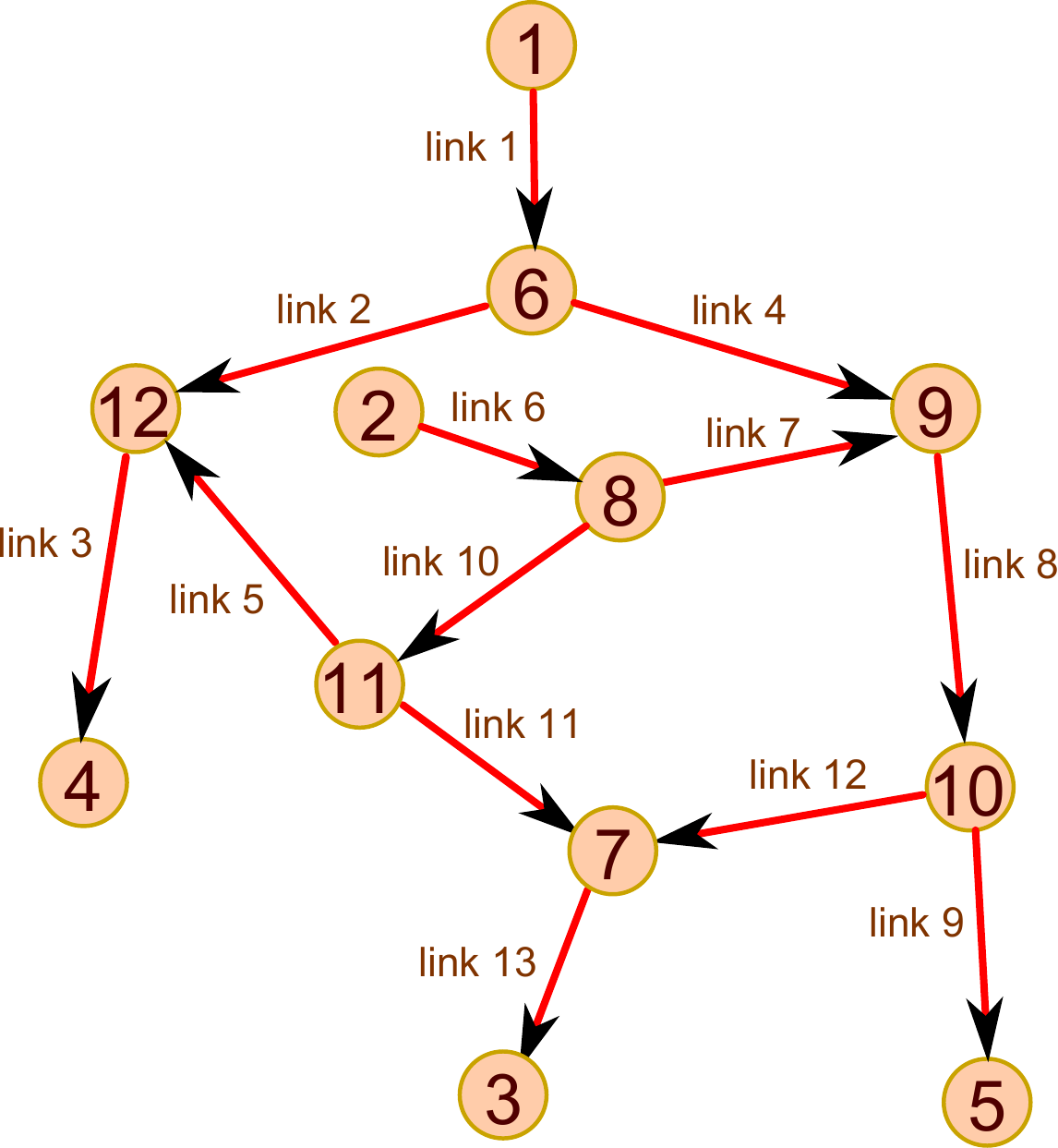}
  \caption{Network for simulation experiment}
  \label{fig:simNet}
\end{figure}

\begin{figure*}[ht!]
  \centerline{
      \hspace{-3em}
    \subfigure[Link 1]{
    \includegraphics[scale = 0.5]{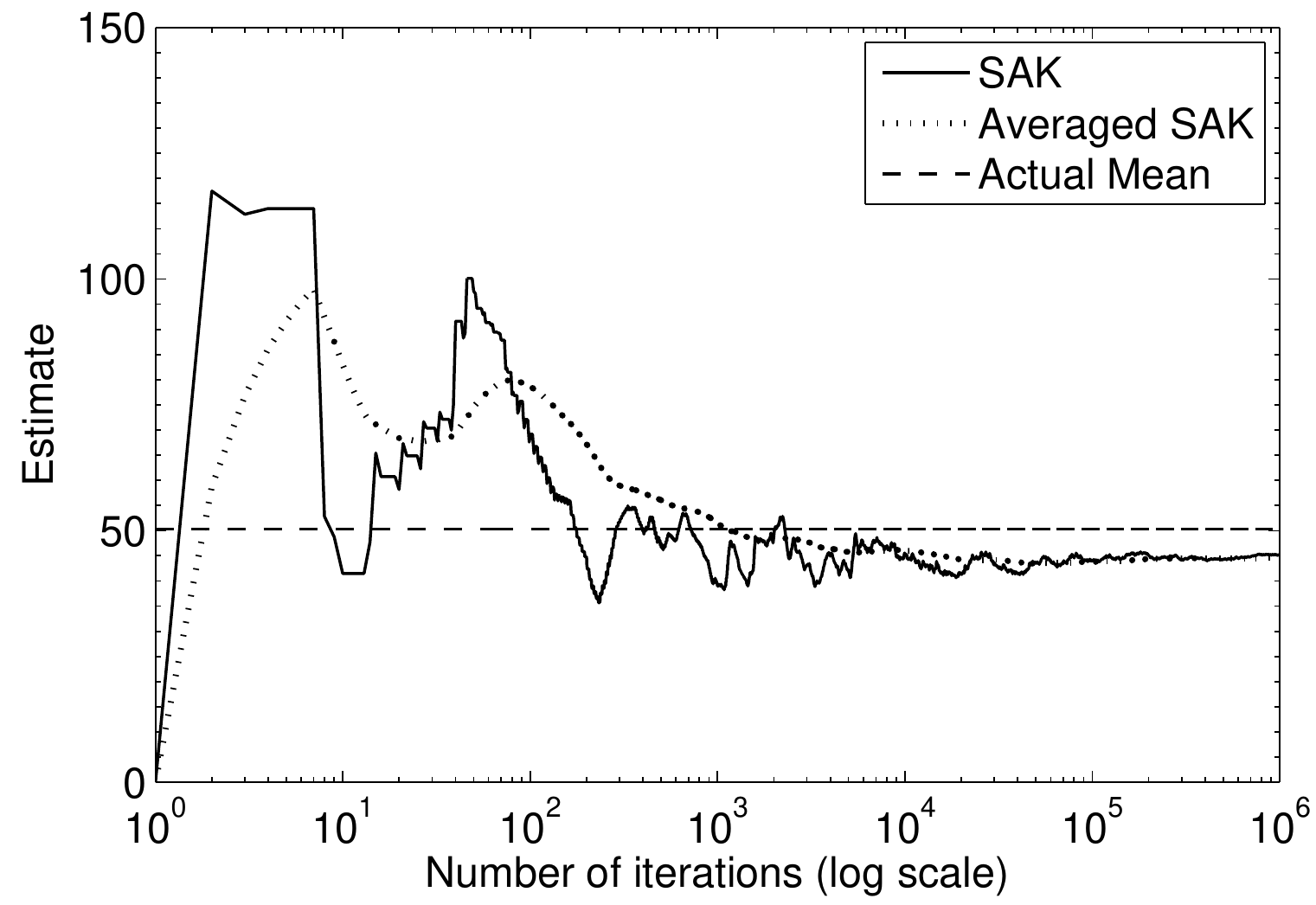}
      \label{fig:link1}
    }
    \hspace{2em}
    \subfigure[Link 3]{
      \includegraphics[scale = 0.5]{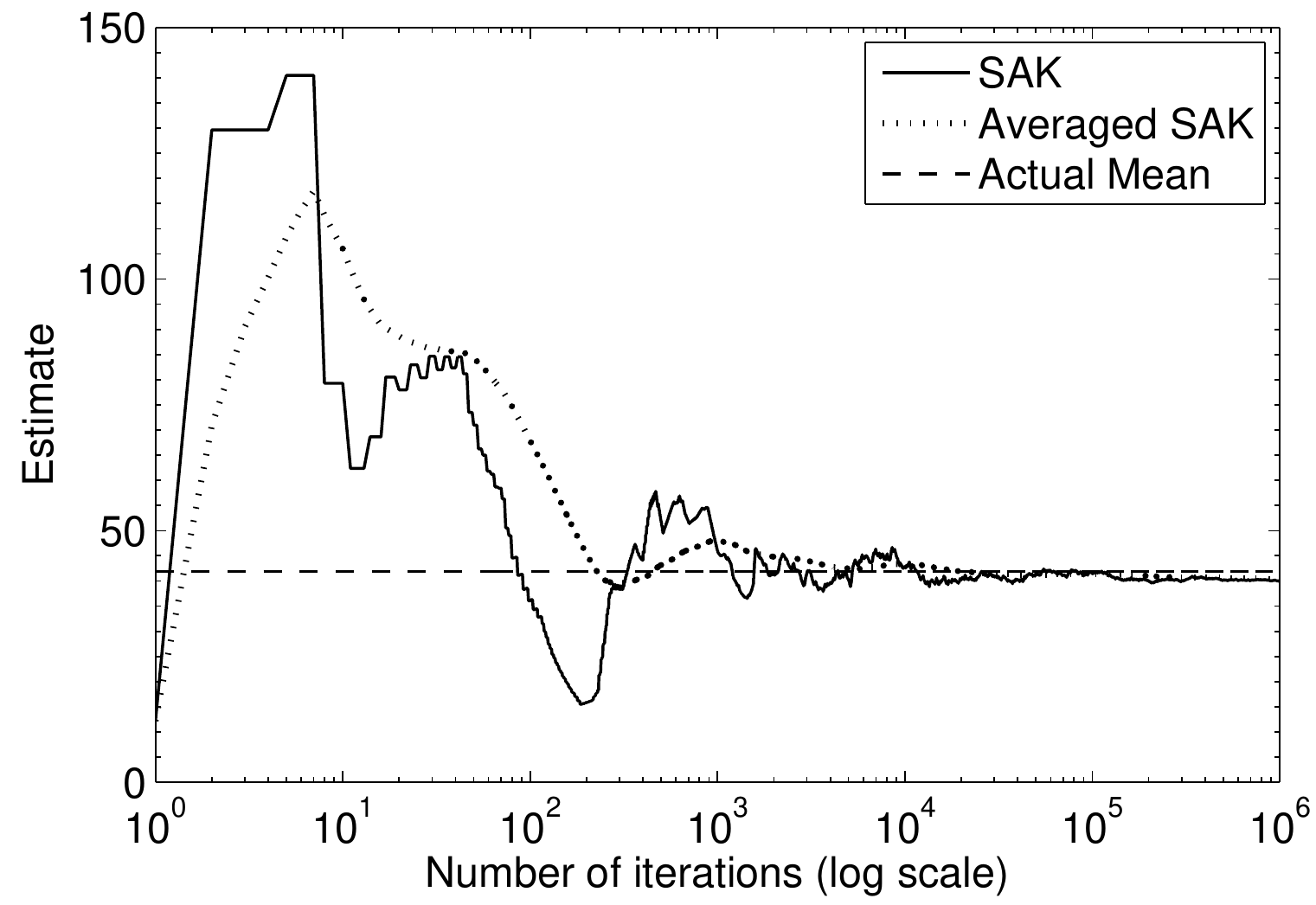}
      \label{fig:link4}
    }
  }
  \caption{Online estimation of expected delay across candidate links 1 and 3 using SAK and Averaged SAK algorithms.}
  \label{fig:Progress}
\end{figure*}

We illustrate the application of SAK algorithm in real time delay tomography for the network of Figure~\ref{fig:simNet}. The goal here is to use the measurements of end-to-end delay experienced by probe packets while traversing different paths in the network to obtain, in real time, the estimates of link delay statistics. 

In the framework of Section~\ref{sec:model}, the experimental setup is as follows. A priori we choose six paths in the network. This is described by the path-link matrix (rows $\mathrel{\widehat{=}}$ paths, columns $\mathrel{\widehat{=}}$ links)

\[
A=
\left[\begin{array}{ccccccccccccc}
1 & 1 & 1 & 0 & 0 & 0 & 0 & 0 & 0 & 0 & 0 & 0 & 0\\
1 & 0 & 0 & 1 & 0 & 0 & 0 & 1 & 1 & 0 & 0 & 0 & 0\\
1 & 0 & 0 & 1 & 0 & 0 & 0 & 1 & 0 & 0 & 0 & 1 & 1\\
0 & 0 & 1 & 0 & 1 & 1 & 0 & 0 & 0 & 1 & 0 & 0 & 0\\
0 & 0 & 0 & 0 & 0 & 1 & 1 & 1 & 1 & 0 & 0 & 0 & 0\\
0 & 0 & 0 & 0 & 0 & 1 & 0 & 0 & 0 & 1 & 1 & 0 & 1
\end{array}\right].
\]

Its entry $a_{ij}$ is one if link $j$ is present on path $i.$ Thus, row four denotes the path that connects the nodes 2-8-11-12-4. The delay a probe packet experiences while traversing link $j$ is a random variable $X(j)$ with arbitrary non-negative distribution. The delay across path $i$ is $Y(i) = \langle a_i, X \rangle + W(i),$ where $W(1), \ldots, W(6)$ are IID standard Gaussian random variables denoting measurement error. We generate a million probe packets, where the $k^{th}$ packet is sent along a path whose index, denoted $Z_k,$ is chosen uniformly randomly from $\{1, \ldots, 6\}.$ Thus each path gets about 167,000 samples. We use $\mathcal{Y}_k$ to record the delay, packet $k$ experiences while traversing the path $Z_k.$ We also run our SAK Algorithm of (\ref{eq:SAK}) for a million iterations, first for (\ref{eq:LinearRel}) and then for (\ref{eq:HighMomMatRel}) with $q = 2.$ The chosen start point, the actual value and final estimated value of moments are given in Tables~\ref{tab:simResults1} and \ref{tab:simResults2}. In both cases the initial point satisfies the assumption of Lemma~\ref{lem:QualityInitPoint} and hence the final estimates are close to the actual values. In Table~\ref{tab:simResults2}  we give only a subset of results. Figure~\ref{fig:Progress} compares the real-time estimates of expected delay for candidate links $1$ and $3$ obtained using the SAK algorithm and the averaged SAK algorithm. The iterates of the averaged SAK algorithm are samples averages of the SAK algorithm iterates. Observe that, although we run the simulation for a million packets, the estimates are very nearly the true values in after about 300 iterations. Also, note that the error in the estimates does not decrease monotonically. This is because of the direct use of noisy measurements. The fluctuations, however, get suppressed as the stepsizes decrease with iterations.

\begin{table}[ht!]
  \caption{Actual, Final estimated value of Expected Link Delay}
  \label{tab:simResults1}
  \vspace{0.5ex}
  \centering
  \begin{tabular}{|cccc|}
    \hline
    \multirow{2}{*}{LinkId} & Initial & True expected & Final \vspace{-1.5ex}\\
     & guess & delay & estimate\tabularnewline
    \hline
    \hline
    1 & 00.00  & 50.25  & 45.09\tabularnewline
    \hline
    2 & 00.00  & 26.32  & 33.17\tabularnewline
    \hline
    3 & 12.15  & 41.84  & 39.96\tabularnewline
    \hline
    4 & 00.00  & 09.10  & 11.92\tabularnewline
    \hline
    5 & 25.34  & 23.04  & 19.98\tabularnewline
    \hline
    6 & 00.00  & 48.08  & 46.87\tabularnewline
    \hline
    7 & 00.00  & 41.49  & 39.05\tabularnewline
    \hline
    8 & 00.00  & 49.75  & 50.97\tabularnewline
    \hline
    9 & 00.00  & 34.72  & 37.34\tabularnewline
    \hline
    10 & 00.00  & 03.78  &  07.82\tabularnewline
    \hline
    11 & 28.86  &44.05   & 42.06\tabularnewline
    \hline
    12 & 39.90  & 48.54  & 53.54\tabularnewline
    \hline
    13 & 00.00  & 29.07  & 26.82 \tabularnewline
    \hline
  \end{tabular}
\end{table}
\vspace{1ex}

\begin{table}[ht]
  \caption{Actual, Final estimated value of some $2^{nd}$ Order Moments}
  \label{tab:simResults2}
  \vspace{1ex}
  \centering
  \begin{tabular}{|cccc|}
    \hline
    \multirow{2}{*}{Moment} & Initial & True  & Final \vspace{-1.5ex}\\
     & guess & value & estimate\tabularnewline
    \hline
    \hline
    \mbox{$\mathbb{E}(X_1^2)$}     & 17388  & 20539  & 20570\tabularnewline
    \hline
    \mbox{$\mathbb{E}(X_4^2)$}     & 0      & 277.85 & 286.29\tabularnewline
    \hline
    \mbox{$\mathbb{E}(X_3X_{10})$} & 15985  & 158.83 & 164.34\tabularnewline
    \hline
	\mbox{$\mathbb{E}(X_8X_{12})$} & -126   & 2427.8 & 2390.5\tabularnewline
    \hline
  \end{tabular}
\end{table}

\bibliographystyle{plain}
\bibliography{tomography}

\begin{thebibliography}{10}

\bibitem{Abramsson98}
T.~Abrahamsson.
\newblock Estimation of origin-destination matrices using traffic counts--a
  literature survey.
\newblock {\em IIASA Interim Report IR-98-021/May}, 27:76, 1998.

\bibitem{Adams00}
A.~Adams, T.~Bu, T.~Friedman, J.~Horowitz, D.~Towsley, R.~Caceres, N.~Duffield,
  F.~Presti, S.~B. Moon, and V.~Paxson.
\newblock The use of end-to-end multicast measurements for characterizing
  internal network behavior.
\newblock {\em Communications Magazine, IEEE}, 38(5):152--159, 2000.

\bibitem{Ananth12}
V.~Anantharam and V.~S. Borkar.
\newblock Stochastic approximation with long range dependent and heavy tailed
  noise.
\newblock {\em Queueing Systems}, 71(1-2):221--242, 2012.

\bibitem{Astrom83}
K.J. {\AA}str{\"o}m.
\newblock Theory and applications of adaptive control--a survey.
\newblock {\em Automatica}, 19(5):471--486, 1983.

\bibitem{Bell91}
M.~G.~H. Bell.
\newblock The estimation of origin-destination matrices by constrained
  generalised least squares.
\newblock {\em Transportation Research Part B: Methodological}, 25(1):13--22,
  1991.

\bibitem{Benveniste90}
A.~Benveniste, M.~M{\'e}tivier, and P.~Priouret.
\newblock {\em Adaptive algorithms and stochastic approximations}.
\newblock Springer Publishing Company, Incorporated, 1990.

\bibitem{Borkar08}
V.~S. Borkar.
\newblock {\em Stochastic approximation: a dynamical systems viewpoint}.
\newblock Cambridge University Press Cambridge, 2008.

\bibitem{Caceres99}
R.~C{\'a}ceres, N.~Duffield, J.~Horowitz, and D.~Towsley.
\newblock Multicast-based inference of network-internal loss characteristics.
\newblock {\em Information Theory, IEEE Transactions on}, 45(7):2462--2480,
  1999.

\bibitem{Castro04}
R.~Castro, M.~Coates, G.~Liang, R.~Nowak, and B.~Yu.
\newblock Network tomography: recent developments.
\newblock {\em Statistical science}, pages 499--517, 2004.

\bibitem{Chong01}
E.~Chong and H.~Zak.
\newblock {\em An introduction to optimization}.
\newblock John Wiley \& Sons, 2001.

\bibitem{Coates02b}
A.~Coates, A.~Hero III, R.~Nowak, and B.~Yu.
\newblock Internet tomography.
\newblock {\em Signal Processing Magazine, IEEE}, 19(3):47--65, 2002.

\bibitem{Derevitskii74}
D.~Derevitskii and A.~Fradkov.
\newblock Two models for analyzing the dynamics of adaptation algorithms.
\newblock {\em Automation and Remote Control}, 35(1):59--67, 1974.

\bibitem{Downey99}
A.~Downey.
\newblock Using pathchar to estimate internet link characteristics.
\newblock In {\em ACM SIGCOMM Computer Communication Review}, volume 29(4),
  pages 241--250. ACM, 1999.

\bibitem{Duffield01}
N.~Duffield, J.~Horowitz, F.~Presti, and D.~Towsley.
\newblock Network delay tomography from end-to-end unicast measurements.
\newblock In {\em Evolutionary Trends of the Internet}, pages 576--595.
  Springer, 2001.

\bibitem{Duffield02}
N.~Duffield, J.~Horowitz, F.~Presti, and D.~Towsley.
\newblock Multicast topology inference from measured end-to-end loss.
\newblock {\em Information Theory, IEEE Transactions on}, 48(1):26--45, 2002.

\bibitem{Feldmann01}
A.~Feldmann, A.~Greenberg, C.~Lund, N.~Reingold, J.~Rexford, and F.~True.
\newblock Deriving traffic demands for operational ip networks: Methodology and
  experience.
\newblock {\em IEEE/ACM Transactions on Networking (ToN)}, 9(3):265--280, 2001.

\bibitem{Freris12}
N.~Freris and A.~Zouzias.
\newblock Fast distributed smoothing of relative measurements.
\newblock In {\em Decision and Control (CDC), 2012 IEEE 51st Annual Conference
  on}, pages 1411--1416. IEEE, 2012.

\bibitem{Galantai04}
A.~Gal{\'a}ntai.
\newblock {\em Projectors and projection methods}, volume~6.
\newblock Springer, 2004.

\bibitem{Robbins51}
Robbins H and S.~Monro.
\newblock A stochastic approximation method.
\newblock {\em The Annals of Mathematical Statistics}, pages 400--407, 1951.

\bibitem{Kaczmarz37}
S.~Kaczmarz.
\newblock Angen{\"a}herte aufl{\"o}sung von systemen linearer gleichungen.
\newblock {\em Bulletin International de l’Academie Polonaise des Sciences et
  des Lettres}, 35:355--357, 1937.

\bibitem{Kushner03}
H.~Kushner and G.~Yin.
\newblock {\em Stochastic approximation algorithms and applications}.
\newblock Springer New York, 2003.

\bibitem{Lawrence07}
E.~Lawrence, G.~Michailidis, and V.~Nair.
\newblock Statistical inverse problems in active network tomography.
\newblock {\em Lecture Notes-Monograph Series}, pages 24--44, 2007.

\bibitem{Leventhal10}
D.~Leventhal and A.~Lewis.
\newblock Randomized methods for linear constraints: convergence rates and
  conditioning.
\newblock {\em Mathematics of Operations Research}, 35(3):641--654, 2010.

\bibitem{Ljung77}
L.~Ljung.
\newblock Analysis of recursive stochastic algorithms.
\newblock {\em Automatic Control, IEEE Transactions on}, 22(4):551--575, 1977.

\bibitem{Maher83}
M.~Maher.
\newblock Inferences on trip matrices from observations on link volumes: a
  bayesian statistical approach.
\newblock {\em Transportation Research Part B: Methodological}, 17(6):435--447,
  1983.

\bibitem{Moustakides98}
G.~Moustakides.
\newblock Exponential convergence of products of random matrices: Application
  to adaptive algorithms.
\newblock {\em International Journal of Adaptive Control and Signal
  Processing}, 12(7):579--597, 1998.

\bibitem{Parks92}
P.~Parks and J.~Militzer.
\newblock A comparison of five algorithms for the training of cmac memories for
  learning control systems.
\newblock {\em Automatica}, 28(5):1027--1035, 1992.

\bibitem{Polyak92}
B.~Polyak and A.~Juditsky.
\newblock Acceleration of stochastic approximation by averaging.
\newblock {\em SIAM Journal on Control and Optimization}, 30(4):838--855, 1992.

\bibitem{Popa04}
C.~Popa and R.~Zdunek.
\newblock Kaczmarz extended algorithm for tomographic image reconstruction from
  limited-data.
\newblock {\em Mathematics and Computers in Simulation}, 65(6):579--598, 2004.

\bibitem{Prasad03}
R.~Prasad, C.~Dovrolis, M.~Murray, and K.~Claffy.
\newblock Bandwidth estimation: metrics, measurement techniques, and tools.
\newblock {\em Network, IEEE}, 17(6):27--35, 2003.

\bibitem{LoPresti02}
F.~Presti, N.~Duffield, J.~Horowitz, and D.~Towsley.
\newblock Multicast-based inference of network-internal delay distributions.
\newblock {\em Networking, IEEE/ACM Transactions on}, 10(6):761--775, 2002.

\bibitem{Richalet78}
J.~Richalet, A.~Rault, J.~Testud, and J.~Papon.
\newblock Model predictive heuristic control: Applications to industrial
  processes.
\newblock {\em Automatica}, 14(5):413--428, 1978.

\bibitem{Sherali94}
H.~Sherali, R.~Sivanandan, and A.~Hobeika.
\newblock A linear programming approach for synthesizing origin-destination
  trip tables from link traffic volumes.
\newblock {\em Transportation Research Part B: Methodological}, 28(3):213--233,
  1994.

\bibitem{Strohmer09}
T.~Strohmer and R.~Vershynin.
\newblock A randomized {Kaczmarz} algorithm with exponential convergence.
\newblock {\em Journal of Fourier Analysis and Applications}, 15(2):262--278,
  2009.

\bibitem{Vardi96}
Y.~Vardi.
\newblock Network tomography: Estimating source-destination traffic intensities
  from link data.
\newblock {\em Journal of the American Statistical Association},
  91(433):365--377, 1996.

\bibitem{Zhang03a}
Y.~Zhang, M.~Roughan, N.~Duffield, and A.~Greenberg.
\newblock Fast accurate computation of large-scale ip traffic matrices from
  link loads.
\newblock In {\em ACM SIGMETRICS Performance Evaluation Review}, volume 31(1),
  pages 206--217. ACM, 2003.

\bibitem{Zhang03b}
Y.~Zhang, M.~Roughan, C.~Lund, and D.~Donoho.
\newblock An information-theoretic approach to traffic matrix estimation.
\newblock In {\em Proceedings of the 2003 conference on Applications,
  technologies, architectures, and protocols for computer communications},
  pages 301--312. ACM, 2003.

\bibitem{Zouzias12}
A.~Zouzias and N.~Freris.
\newblock Randomized extended {Kaczmarz} for solving least-squares.
\newblock {\em arXiv preprint arXiv:1205.5770}, 2012.

\end{thebibliography}

\end{document}